\documentclass{amsart}

\usepackage{amsmath}
\usepackage{amsthm}
\usepackage{stmaryrd}
\usepackage{amsfonts,mathrsfs}
\title[Quotients of squared norms]{Geometry of Hermitian Algebraic Functions.  \\
Quotients of squared norms}
\author{Dror Varolin}
\thanks{Partially supported by NSF grant DMS-0400909}
\address{Department of Mathematics\\ Stony Brook University\\ Stony Brook, NY 11794}
\email{dror@math.sunysb.edu}

\newcommand{\noi}{\noindent}

\newcommand{\ms}{\medskip}

\newcommand{\gmod}[1]{{\talloblong #1 \talloblong}}

\newcommand{\ch}{{\mathcal H}}

\newcommand{\co}{{\mathcal O}}

\newcommand{\cq}{{\mathcal Q}}

\newcommand{\cz}{{\mathcal Z}}

\newcommand{\sa}{{\mathscr A}}
\newcommand{\sB}{{\mathscr B}}
\newcommand{\sC}{{\mathscr C}}

\newcommand{\sP}{{\mathscr P}}

\newcommand{\vp}{\varphi}

\newcommand{\C}{{\mathbb C}}
\newcommand{\h}{{\mathbb H}}
\newcommand{\K}{{\mathbb K}}
\newcommand{\N}{{\mathbb N}}
\newcommand{\p}{{\mathbb P}}

\newcommand{\Z}{{\mathbb Z}}

\newcommand{\di}{\partial}
\newcommand{\dbar}{\bar \partial}

\newcommand{\ii}{\sqrt{-1}}

\newcommand{\emb}{\hookrightarrow}
\newcommand{\tensor}{\otimes}

\begin{document}
\maketitle

\theoremstyle{plain}
\newtheorem{thm}{\sc Theorem}
\newtheorem{lem}{\sc Lemma}[section]
\newtheorem{o-thm}[lem]{\sc Theorem}
\newtheorem{prop}[lem]{\sc Proposition}
\newtheorem{cor}[lem]{\sc Corollary}

\theoremstyle{definition}
\newtheorem{conj}[lem]{\sc Conjecture}
\newtheorem{defn}[lem]{\sc Definition}
\newtheorem{qn}[lem]{\sc Question}
\newtheorem{ex}[lem]{\sc Example}

\theoremstyle{remark}
\newtheorem*{rmk}{\sc Remark}
\newtheorem*{ack}{\sc Acknowledgment}

\setcounter{section}{-1}

\section{Introduction}\label{intro}

In this paper we characterize those Hermitian algebraic functions that are quotients of squared norms of holomorphic mappings.  Our approach is to use the resolution of singularities, followed by basic algebraic geometry of big line bundles, to reduce the problem to elementary properties of the Bergman kernel for spaces of sections of large tensor powers of a positive holomorphic line bundle.

Recall that if $M$ is a complex manifold and $M^{\dagger}$ is the complex conjugate manifold, then a function $\ch : M \times M^{\dagger} \to \C$ is said to be {\it Hermitian} if $\ch$ is holomorphic and $\ch (z,\bar w) = \overline{\ch (w,\bar z)}$.   By polarization, any Hermitian function is determined by its (real) values $\ch (z,\bar z)$ along the ``diagonal".  If $\ch (z,\bar z) \ge 0$ for all $z \in M$, we shall say that $\ch$ is {\it non-negative}.

If the manifold $M$ in question is the total space of the dual of a holomorphic line bundle $F \to X$ on a projective manifold $X$, then a Hermitian function $P :F^* \times (F^*)^{\dagger} \to \C$ is called {\it Hermitian algebraic} if for any $\lambda \in \C$ and $v,w \in F^*$, $P(\lambda \cdot v,\bar w) = \lambda P(v, \bar w)$. (Scalar multiplication on $F^*$ is along fibers.) Equivalently, a Hermitian algebraic function $P$ is a function of the form 
\begin{equation}\label{herm-rep}
P(v,\bar w) = \sC _{\alpha \bar \beta} \left < s^{\alpha} ,v\right > \overline{\left < s^{\beta} ,w\right >},
\end{equation}
where $(\sC _{\alpha \bar \beta})$ is a Hermitian matrix and $\{s^{\alpha}\}$  are global holomorphic sections of $F \to X$.  (Here and in the rest of the paper, the usual summation convention of summing over repeated indices is in force.)  We shall employ the notation
\[
P= \sC _{\alpha \bar \beta} s^{\alpha } \bar s ^{\beta}.
\]

\begin{ex}
If $\h \to \p _n$ is the line bundle associated to a hyperplane
section and $d$ is a positive integer, then the global sections of $\h ^{\tensor d}$ are homogeneous polynomials of degree $d$ in the
homogeneous coordinates $z=[z^0,...,z^n]$.  In this setting a
Hermitian algebraic function is a bihomogeneous polynomial
\[
P(z,\bar z) = \sum _{|\alpha|=|\beta|=d} C_{\alpha \bar \beta} z^{\alpha} \bar z ^{\beta},
\]
where $C_{\alpha \bar \beta} = \overline{C _{\beta \bar \alpha}}$.
\end{ex}

Upon diagonalizing the Hermitian matrix $\sC _{\alpha \bar \beta}$, we find that there exist independent sections $f^1,...,f^k, g^1,...,g^{\ell}$ such that 
\begin{equation}\label{dos}
P = \sum _{i=1} ^k f^i \bar f^i - \sum _{j=1}^{\ell} g^j \bar g^j.
\end{equation}
We define the {\it modulus} of $P$ to be the Hermitian algebraic function
\[
\gmod{P} = \sum _{i=1} ^k f^i \bar f^i + \sum _{j=1}^{\ell} g^j \bar g^j.
\]
We warn the reader that $\gmod{P}$ is not canonically determined by $P$, but also depends on the basis $\{ f^1,..., f^k, g^1,..., g^{\ell} \}$.

A Hermitian algebraic function $P :F^* \times (F^*)^{\dagger} \to \C$ is called a {\it sum of squared norms} if $P$ can be represented in the form (\ref{herm-rep}) with $(\sC _{\alpha \bar \beta})$ positive semi-definite.  The decomposition \eqref{dos} then reduces to 
\[
P(v,\bar v) = \sum _{\alpha =1} ^k |\left < f^{\alpha} , v\right >
|^2.
\]

If $P$ is a Hermitian algebraic function and there is a sum of
squared norms $Q$ such that the Hermitian algebraic function $R = QP$ is a sum of squared norms, then we say that $P$ is a {\it quotient of squared norms}.

Clearly a quotient of squared norms is non-negative.  A Hermitian analogue of Hilbert's 17th problem, posed by D'Angelo \cite{d-carrus,d2}, is to characterize those non-negative Hermitian algebraic functions that are quotients of squared norms.  Unlike the ``real" 17th problem of Hilbert, there are many examples of non-negative Hermitian algebraic functions that are not quotients of squared norms.   The simplest examples are constructed based on the fact that if a Hermitian algebraic function $P$ is a quotient of squared norms, then the zero set of the function $z \mapsto P(z,\bar z)$ must be an analytic set.  
\begin{ex}
In the case $X= \p _1$ and $F = \h ^{\tensor 2}$, the Hermitian algebraic function $P = (|z_0|^2 -|z_1|^2)^2$ is not a quotient of squares, because the set of $x \in \p _1$ such that $P(v,\bar v) = 0$ for some (and thus all) $v \in (\h ^{\tensor 2})^*_x$ is a circle.
\end{ex}
In fact the phenomenon is more complex than this; in Section \ref{qos-section} we mention two examples, the second of which is due to D'Angelo, of non-negative Hermitian algebraic functions that are not quotients of squared norms, but whose zero sets are complex analytic sets.

Quillen \cite{q} and independently Catlin and D'Angelo \cite{cd1} showed that if a Hermitian bihomogeneous polynomial is positive on the unit sphere, then it is a quotient of squared norms.  Later Catlin and D'Angelo \cite{cd3} generalized this result to show that any Hermitian algebraic function that is strictly positive away from the zero section is a quotient of squared norms (Theorem \ref{cd-qos-thm} below).

Catlin and D'Angelo deduce Theorem \ref{cd-qos-thm} from a more powerful theorem (stated as Theorem \ref{cd3-thm} below), which is the main result in the same paper \cite{cd3}.  A key tool they use is a theorem of Catlin \cite{c} regarding the asymptotic expansion of the Bergman kernel.  Zelditch \cite{z} proved an analogous result for the Szeg\" o kernel.  (Both Catlin and Zelditch used their results to settle a conjecture of Tian regarding approximation of K\" ahler-Einstein metrics.)  The Catlin-Zelditch theorem is a version of the well-known generalization, due to Boutet de Monvel and Sj\" ostrand \cite{bs}, of the celebrated theorem of Fefferman \cite{f}
on the asymptotic expansion of the Bergman kernel.

In this paper, we completely characterize all those non-negative Hermitian symmetric functions that are quotients of squared norms.  The following is our main result.

\begin{thm}\label{main}
A non-negative Hermitian algebraic function $P$ is a quotient of squared norms if and only if there exists a constant $C$ such that 
\[
\text{for all } v \in F^* - \{ v \ ;\ P(v,\bar v) = 0\}, \quad \frac{\gmod{P}(v,\bar v)}{P(v,\bar v)} \le C.
\]
\end{thm}

\begin{rmk}
As we said, $\gmod{P}$ is not uniquely defined by $P$.  However, the necessary and sufficient condition in Theorem \ref{main} is canonically determined by $P$ (Proposition \ref{canon}).
\end{rmk}

\begin{rmk}
D'Angelo has also characterized quotients of squared norms in
\cite{d2}, at least in the case of bihomogeneous Hermitian
polynomials.  For the case of bihomogeneous polynomials on $\C ^2$ (which corresponds here to the case $X= \p _1$) his result is easily seen to be equivalent to ours, but in higher dimensions the nature of his necessary and sufficient condition is very different from ours.  
\end{rmk}

Our approach to the proof of Theorem \ref{main} is to apply Hironaka's Theorem on the resolution of singularities, along with some methods from elementary algebraic geometry, to reduce the problem to the case in which the aforementioned result of Catlin-D'Angelo, namely Theorem \ref{cd3-thm} below, applies.

Because of its role in our proof, for the sake of self containment we give a short proof of Theorem \ref{cd3-thm} in the final section of this paper.  We claim essentially no originality here; the main idea of the proof we give,  essentially the same idea of Catlin-D'Angelo, is to make use of the asymptotic expansion for the Bergman kernel.  However, instead of using the theorem of Catlin about compact perturbations of the Bergman kernel of the disk bundle $\{ R(v,\bar v) < 1\}$, we use directly the asymptotic expansion for the Bergman kernel associated to $L^2$ spaces of sections of large powers of holomorphic line bundles.  (The latter, while a consequence of the former, has recently been established directly by Berman-Berndtsson-Sj\" ostrand \cite{bbs} using elementary methods.)  This allows us to give a more direct argument, and avoid compact integral operators.  In fact, the last part of the proof is instead closer in spirit to Quillen's proof \cite{q} of the special case of Hermitian polynomials.

\tableofcontents

\begin{ack}
I am grateful to John D'Angelo for introducing me to the study of Hermitian algebraic functions and posing to me the general problem of characterizing quotients of squared norms, and for his generous help with the preparation of this manuscript.  I am also grateful to Dave Catlin for pointing out an error in an earlier version of the paper, and to Mark de Cataldo, Rob Lazarsfeld, Laszlo Lempert and Jeff McNeal for stimulating discussions.  Finally, I'm grateful to the referee for finding some errors in an earlier version of the manuscript, and for suggestions that I hope have vastly improved the presentation.
\end{ack}

\section{Some notation and background from complex and algebraic
geometry}\label{background-section}

In this section we discuss some preliminary ideas needed in this work.   We also take the opportunity to establish notation that will be used in the rest of the paper.

\medskip

\noi {\bf Holomorphic  line bundles.}
Let $X$ be a complex manifold and $E \to X$ a holomorphic line bundle.  We denote by $H^0(X,E)$ the space of global holomorphic sections of $E$.  By the compactness of $X$, $H^0(X,E)$ is a finite-dimensional vector space.

Recall that a $\Z$-divisor (or simply divisor) on $X$ is a $\Z$-linear combination of irreducible hypersurfaces on $X$.  Given a divisor $D$ on $X$, there is a holomorphic line bundle $L_D \to X$ and a meromorphic section $s_D :X \to L_D$ of $L_D$ such that the divisor of $s_D$ (i.e. the zero divisor minus the polar divisor) is precisely $D$.  In particular, if $D$ is effective then $s_D$ is holomorphic.  The line bundle $L_D$ is defined by the transition functions
\[
\frac{f_{\alpha}}{f_{\beta}} \quad \text{on }U_{\alpha} \cap U_{\beta},
\]
where $f_{\alpha}$ is the local defining function for $D$ on
$U_{\alpha}$.  The section $s_D$ is defined over $U_{\alpha}$ by the local representative $f_{\alpha}$.  

\medskip

\noi {\bf Hermitian metrics.}  In the last section of the paper, when we consider $L^2$ spaces of sections of a holomorphic line bundle, we shall need the notion of a Hermitian metric for a holomorphic line bundle $E \to X$.   Recall that a Hermitian metric $h$ on a complex line bundle $E \to X$ is a smooth family of Hermitian metrics on the fibers of $E$.  Let $U$ be an open set in $X$ and assume there is a nowhere zero section $e \in H^0(U, E|U)$.  Given a local section $v$ of $E$ over the neighborhood $U$, we have $v = fe$ for some function $f$, and thus 
\[
h(v,v) = |f|^2 e^{-\vp_U}, \quad \text{where} \quad \vp _U = - \log h(e,e).
\]
If we do this for an open cover $\{ U_j\}$, the $(1,1)$-forms $\ii \di \dbar \vp _{U_j}$ agree on overlaps on their domains of definition, and thus define a global form on $X$ called the {\it curvature} of $h$.  It is common to omit the subscript $U$ and employ the abusive notation 
\[
h = e^{-\vp} \quad \text{and} \quad h(v,v) = |v|^2 e^{-\vp}.
\]
The curvature form is then denoted $\ii \di \dbar \vp$.

Although we will not use them in a fundamental way, it is convenient to make use of the notion of singular Hermitian metric for a holomorphic line bundle.  A {\it singular} Hermitian metric is similar to a Hermitian metric except that the functions $\vp$ are required to be only locally integrable.  The curvature
\[
\ii \di \dbar \vp
\]
of a singular metric is thus a well-defined $(1,1)$-current.  We will say
that the singular metric has positive (resp. non-negative)
curvature if the current $\ii \di \dbar \vp$ is positive definite
(resp. positive semi-definite) in the sense of currents, i.e., if
the metric is given locally by $e^{-\vp}$ with $\vp$ strictly
plurisubharmonic (resp. plurisubharmonic).

\begin{rmk}
With this unfortunately common notation, a smooth Hermitian metric is also a singular Hermitian metric.  When emphasis is required, we will add the adjective {\it smooth}.
\end{rmk}

\medskip

\noi {\bf Resolution of singularities.}
We shall need Hironaka's Theorem on log resolution of singularities of linear systems \cite{h}.  The precise result we need is the following.

\begin{o-thm}\label{log-res}{\rm (Hironaka)}
Let $F \to X$ be a holomorphic line bundle and $V \subset
H^0(X,F)$ a non-zero finite-dimensional subspace.  Then there is a birational map $\mu : \tilde X \to X$ with $\tilde X$ smooth, an effective normal crossing divisor $\sum D_j$, and non-negative integers $m_j, e_j$ such that, with $D = \sum m_j D_j$ and $s_D$ the canonical holomorphic section of the line bundle associated to $D$ whose zero divisor is $D$,
\begin{enumerate}
\item[(i)] $J := {\rm Exc}(\mu) = \sum e_j D_j$, 

\item[(ii)] $\mu ^* V = s_D \tensor W$, and 

\item[(iii)] $W \subset H^0(\tilde X, \mu ^* F - D) $ is base point free.
\end{enumerate}
\end{o-thm}

\noi Recall that a subspace $V \subset H^0(X,L)$ of global sections of a holomorphic line bundle is base point free (or simply, free) if for each $x \in X$ there is a section $s \in V$ such that $s(x) \neq 0$, and that $\text{Exc}(\mu)$ is the smallest divisor outside of which the map $\mu$ is an isomorphism.  The map $\mu$ in Hironaka's Theorem is called a log resolution of $V$.

\begin{rmk}
Hironaka's original theorem actually says more than this, providing $\mu$ as a composition of blow-ups along smooth centers.  However, the above theorem has a significantly more elementary proof.  See, for example, \cite{bp,ad,p} or the discussion in \cite[Page I-241]{rlaz}.
\end{rmk}

\medskip

\noi {\bf Big line bundles.}
The notion of ample line bundle is not preserved under birational maps.  However, a certain kind of positivity does persist under birational transformation, and this leads to the notion of a big line bundle.

\begin{defn}
Let $A \to X$ be a holomorphic line bundle over a projective
algebraic manifold of dimension $n$. We say that $A$ is {\emph big} if for all sufficiently large integers $m$ the
map
\[
\phi _{A^{\tensor m}} : X \dashrightarrow \p ( H^0 (X,A ^{\tensor m}) )^*
\]
sending $p$ to the hyperplane
\[
\phi _{A^{\tensor m}} (p) := \{ s \in H^0(X,A^{\tensor m})\ | \ s(p) = 0 \}
\]
is birational.
\end{defn}

It is not hard to see that the notion of big line bundle is preserved under birational transformation.  It is a consequence of Kodaira's embedding theorem that every ample line bundle is big.  For us, the most relevant characterization of big is the following fact, due to Kodaira.  For a proof, see \cite[Corollary 2.2.7]{rlaz}.

\begin{prop}\label{big-prop}
If $A \to X$ is big, then there is an integer $N>0$ and an ample line bundle $L \to X$ such that line bundle $Z=A ^{\tensor N} \tensor (L^*)$ is effective, i.e., it has a not-identically zero holomorphic section.
\end{prop}

\begin{rmk}
By taking a large tensor power of the big line bundle $A$, we may assume the line bundle $L$ in Proposition \ref{big-prop} is very ample.
\end{rmk}

\section{Elementary Theory of Hermitian algebraic functions}\label{g-metrics-section}

We now give a more thorough treatment of the definitions and elementary facts regarding Hermitian algebraic functions.

\medskip

\noi {\bf Hermitian algebraic functions and their support spaces.}
The following definition was already given.

\begin{defn}
Let $F \to X$ be a holomorphic line bundle.  A Hermitian algebraic function is a function $P : F^* \times (F^*)^{\dagger} \to \C$ of the form 
\[
P(v,\bar w) = \sC _{\alpha \bar \beta} \left < s^{\alpha} , v\right > \overline{\left < s^{\beta},w\right > },
\]
where $\{ s^{\alpha}\}$ is a basis of $H^0(X,F)$ and $\sC _{\alpha \bar \beta}$ a Hermitian matrix on $\C ^N$, $N = \dim _{\C} H^0(X,F)$.  The collection of Hermitian algebraic functions associated to the line bundle $F$ is denoted
\[
\sP _0 (X,F).
\]
\end{defn}

\begin{rmk}
By polarization, the function $P$ is completely determined by its values along the set 
\[
\{ (v,\bar v)\ ;\ v \in F^* \} \subset F^* \times (F^*)^{\dagger}.
\]
We will use the notation 
\[
P= \sC _{\alpha \bar \beta} s^{\alpha} \bar s ^{\beta}.
\]
\end{rmk}

Let $U = U_{\alpha} ^{\beta}$ be a unitary matrix such that 
\[
U^{\dagger}\sC U = (U_{\alpha} ^{\gamma} \sC _{\gamma \bar \nu} (U^{\dagger})^{\bar \nu}_{\bar \beta})  = \left ( 
\begin{array}{ccc|ccc|ccc}
\lambda _1 &0 & 0 & \cdots & 0 & 0 &0&0& 0  \\
0 & \lambda _2 & 0 & \cdots & 0 & 0&0& 0 &0 \\
0 & 0 & \ddots & 0 & 0 & 0&0& 0 &0 \\ \hline
0&0& 0 & - \mu _1 &0 & 0 & \cdots & 0 & 0 \\
0&0& 0 &0 & -\mu _2 & 0 & \cdots & 0 & 0 \\
0&0& 0 &0 & 0 & \ddots & 0 & 0 & 0 \\ \hline
0&0& 0 &0&0& 0 &0&0& 0 \\
0&0& 0 &0&0& 0 &0&\ddots & 0 \\
0&0& 0 &0&0& 0 &0&0& 0 \\
\end{array} \right ), 
\]
where $\lambda _1,..., \lambda _k, \mu _1,..., \mu _{\ell}$ are positive numbers.  (If $\sC$ has non-zero determinant, then there is no final diagonal box of zeros.)  If we define 
\[
f ^1,..., f^k , g^1,..., g^{\ell} \in H^0(X,F)
\]
by
\[
f^i = \sqrt{\lambda _i} U^i _\alpha s^{\alpha}, \quad i=1,..., k \qquad \text{and} \qquad g^j = \sqrt{\mu _j} U^j _\alpha s^{\alpha}, \quad j=1+k,..., \ell +k, 
\]
we find that $f ^1,..., f^k , g^1,..., g^{\ell}$ are linearly independent sections such that 
\[
P = |f^1|^2 + ... + |f^k|^2 - |g^1|^2 - ... - |g^{\ell}|^2.
\]

\begin{defn}
The vector space 
\[
V_P := \text{Span}_{\C} \{ f ^1,..., f^k , g^1,..., g^{\ell}\} \subset H^0(X,F)
\]
is called the {\it support space} of $P$.  A basis $\{ f ^1,..., f^k , g^1,..., g^{\ell} \}$ for $V_P$ such that 
\[
P= |f ^1|^2 +...+|f^k|^2 - |g^1|^2 -...-|g^{\ell}|^2
\]
is called a {\it distinguished basis} for $V_P$.
\end{defn}

\begin{prop}\label{support-space}
The subspace $V_P$ is uniquely determined by the function $P$.
\end{prop}

\begin{proof}
Suppose we have two decompositions 
\[
P= |f ^1|^2 +...+|f^k|^2 - |g^1|^2 -...-|g^{\ell}|^2 = |\hat f ^1|^2 +...+|\hat f^{\hat k}|^2 - |\hat g^1|^2 -...-|\hat g^{\hat \ell}|^2
\]
We can associate to $P$ a unique element of $H^0(X,F) \tensor (H^0(X,F))^{\dagger}$.  Let us complete the basis $f^1,...,f^k ,g^1, ...,g^{\ell}$ to a basis $f^1,...,f^k ,g^1, ...,g^{\ell},h^1 ,...,h ^m$, and consider a dual basis $\xi ^1,..., \xi ^k , \eta ^1 , ..., \eta ^{\ell} , \zeta ^1 ,...,\zeta ^m$ of $H^0(X,F)^*$, i.e., 
\begin{eqnarray*}
&& \left <f^i , \xi ^j \right > = \delta _{ij} , \ \left <g^i , \eta ^j \right > = \delta _{ij} , \  \left <h^i , \zeta ^j \right > = \delta _{ij}, \\
&& \left <f^i , \eta ^j \right > =\left <f^i , \zeta^j \right > =\left <g^i , \xi ^j \right > =\left <g^i , \zeta^j \right > =\left <h^i , \xi ^j \right > =\left <h^i , \eta ^j \right > =0.
\end{eqnarray*}
By pairing the equality
\[
f ^1\tensor \bar f^1 +...+f^k \tensor \bar f^k - g^1\tensor \bar g ^1 -...-g^{\ell}\tensor \bar g^{\ell} = \hat f ^1\tensor \overline{\hat f^1} +...+\hat f^{\hat k} \tensor \overline{\hat f^{\hat k}}  - \hat g^1\tensor \overline{\hat g^1} -...-\hat g^{\hat \ell}\tensor \overline{\hat g^{\hat \ell}}
\]
with the element $\bar \xi ^j$, we see that $f^j$ is a linear combination of the $\hat f^i$' and $\hat g^j$'s.  Similarly, pairing with the element  $\bar \eta ^j$, we see that $g^j$ is also a linear combination of the $\hat f^i$' and $\hat g^j$'s.
\end{proof}

\noi It will be convenient to write $f = (f^1,...,f^k)$ and $g = (g^1,..., g^{\ell})$, as well as
\[
P = ||f||^2 - ||g||^2.
\]

Now suppose we have two distinguished bases $(f,g)$ and $(\hat f , \hat g)$, i.e.,  such that 
\[
P = ||f||^2 - ||g||^2 = ||\hat f||^2 - ||\hat g||^2.  
\]
Let $A$ be an invertible matrix such that 
\[
\binom{f}{g} = A\binom{\hat f}{\hat g}.
\]
Then, with 
\[
\Delta _{k,\ell} =  \left ( 
\begin{array}{cc}
I_k & 0\\ 
0& -I_{\ell}
\end{array}
\right ),
\]
we have 
\begin{equation}\label{uni}
\Delta _{k,\ell} = A^{\dagger} \Delta _{\hat k,\hat \ell} A.
\end{equation}
It follows that the rows of $A$ are mutually orthogonal, and that $k=\hat k$ and $\ell = \hat \ell$.

\begin{defn}\label{gmod-defn}
Let $\{f ^1,..., f^k , g^1,..., g^{\ell}\}$ be a distinguished basis for $V_P$.  Then we define the modulus of $P$ associated to the distinguished basis $\{f ^1,..., f^k , g^1,..., g^{\ell}\}$ to be the Hermitian algebraic function
\[
\gmod{P} := ||f||^2 + ||g||^2.
\]
\end{defn}

Note that $\gmod{P}$ is not uniquely determined by $P$.  Indeed, since
\[
||f||^2 = \frac{\gmod{P} +P}{2} \quad \text{and} \quad ||g||^2 = \frac{ \gmod{P} - P}{2}, 
\]
we see that if two distinguished bases $(f,g)$ and $(\hat f , \hat g)$ yield the same modulus, then $||f||=||\hat f||$ and $||g||=||\hat g||$, and thus $f=U\hat f$ and $g = V\hat g$ for some unitary matrices $U$ and $V$.  However, there are many transformations $A$ satisfying \eqref{uni} which are not block-diagonal with respect to the direct sum decomposition $\C ^{k+\ell}= \C ^k \oplus \C ^{\ell}$.

Nevertheless, the following proposition holds.

\begin{prop}\label{canon}
Let $\{f ^1,..., f^k , g^1,..., g^{\ell}\}$ and $\{\hat f ^1,..., \hat f^k , \hat g^1,..., \hat g^{\ell}\}$ be two distinguished bases for $V_P$.  Then the function 
\[
v \mapsto \frac{||\left < f, v \right > ||^2 + ||\left < g, v\right > ||^2}{||\left < f, v \right > ||^2 - ||\left < g, v\right > ||^2}
\]
is uniformly bounded on $F^* - \{ v\ ;\ P(v,\bar v) =0 \}$ if and only if the same is true of the function 
\[
v \mapsto \frac{||\left < \right . \hspace{-1mm} \hat f, v \hspace{-1mm}\left . \right > ||^2 + ||\left <  \hat g, v\right > ||^2}{||\left <  \right . \hspace{-1mm} \hat f, v \hspace{-1mm}\left . \right > ||^2 - ||\left <  \hat g, v\right > ||^2}
\]
\end{prop}

\begin{proof}
Both statements hold if and only if for every $s \in V_P$, the function
\[
v \mapsto\frac{\left |\left < s , v \right > \right |^2}{P(v,\bar v)}
\]
is uniformly bounded on $F^*  - \{ v\ ;\ P(v,\bar v) = 0\}$.
\end{proof}

Finally, we shall make use of the following notion.

\begin{defn}\label{free-defn}
A G-form $P\in \sP _0(X,F)$ is said to be free if its support space $V_P$ is base point free, i.e., for each $x \in X$ there exists a section $s \in V_P$ such that $s(x) \neq 0$.
\end{defn}

\medskip

\noi {\bf Non-negative Hermitian algebraic functions.}
As mentioned in the Introduction, if the function $v \mapsto P(v,\bar v)$ takes only non-negative values, we call $P$ a non-negative Hermitian algebraic function, and write
\[
P \in \sP _1 (X,F).
\]

A function $P \in \sP _1 (X,F)$ always has zeros, since it is bihomogeneous, and must thus vanish on the zero section $o_{F^*}$ of $F^*$.  Zeros of $P$ lying in the zero section will be considered trivial.

If there is a vector $v \in F^* - o_{F^*}$ such that $P(v,\bar v)=0$, then we will say that $v$ is a non-trivial zero for $P$.  If such a non-trivial zero $v$ lies in the fiber of $F^*$ over some point $x \in X$, then we will write
\[
x \in \cz _P.
\]
The set $\cz _P \subset X$ is called the zero set of $P$.

We let
\[
\sP _1 ^{\sharp} (X,F)
\]
be the set of non-negative Hermitian algebraic functions $P$ whose zero set $\cz _P$ is empty, i.e., $P$ has no non-trivial zeros.  Each $P \in \sP_1^{\sharp} (X,F)$ is said to be {\it positive}.

Eventually we will be interested in non-negative Hermitian algebraic functions with non-trivial zero sets.  An important class of non-negative Hermitian algebraic functions consists of those having only {\it basic zeros.  }

\begin{defn}\label{basic-zeros-defn}
We say that a point $x \in X$ is a basic zero of $P \in \sP _1(X,F)$ if $x \in \cz _P$ and for any $s \in V_P$ the (well defined) function
\begin{eqnarray}\label{bdd}
\frac{|s|^2}{P} : y \mapsto \frac{|\left < s ,v\right >|^2}{P(v,\bar v)}, \quad \pi v = y
\end{eqnarray}
is bounded in the set $U- \cz _P$ for some neighborhood $U$ of $x$ in $X$.  (Here $\pi : F^* \to X$ denotes the projection map.)  We denote by $\sB _P \subset X$ the locus of basic zeros of $P$.  We say that a non-negative Hermitian algebraic function has only basic zeros if for each $s \in V_P$ the function in {\rm (\ref{bdd})} is uniformly bounded on $X - \cz_P$.  
\end{defn}

\begin{rmk}  
\begin{enumerate}

\item Although the set $\sB _P$ of basic zeros is contained in the base locus of $V_P$, these two sets need not agree.

\item Note that $P$ has only basic zeros if and only if with respect to any distinguished basis $\{ f^1,...,f^k, g^1,...,g^{\ell}\}$ of $V_P$, $\gmod{P}/P$ is uniformly bounded on $X-\cz _P$.  Equivalently,  $P$ has only basic zeros if and only if for each distinguished basis $(f,g)$ there is a constant $\lambda \in [0,1)$ such that
\[
||g||^2 \le \lambda ||f||^2.
\]
\end{enumerate}
\end{rmk}

\medskip

\noi {\bf Factoring basic zeros in codimension 1.}  
Let $s$ be a global holomorphic section of a line bundle $F \to X$. Suppose there is a complex hypersurface $D$ on $X$ such that $s |D \equiv 0$.  Thinking of $D$ as a prime divisor, there is a line bundle $L_D$ associated to $D$, and a holomorphic section $s_D$ whose zero locus is precisely $D$.  It follows that $s_D$ divides $s$.  That is to say, the section $s/s_D$ of the line bundle $F\tensor L_D^*$ is holomorphic.  This construction extends to effective divisors $D$ as long as $s$ vanishes on $D$, counting multiplicity.

If all of the sections lying in the support space $V_P$ of a some $P \in \sP _1 (X,F)$ vanish on an effective divisor $D \subset X$, then one can replace $F$ with the line bundle $F\tensor L_D^*$ and $P$ with $\tilde P = P / |s_D|^2$.  By taking $D$ to be the largest divisor in $X$ such that every section of $V_P$ vanishes on $D$, we would obtain in this way a Hermitian algebraic function $\tilde P$ whose basic zero locus has no codimension 1 component.

Thus we will assume from here on out, without loss of generality, that

\begin{center}
{\tt All of our Hermitian algebraic functions have  no basic zeros in codimension 1.}
\end{center}

\medskip

\noi {\bf Associated singular Hermitian metric.}
To every $P \in \sP _1 (X,F)$ we can associate an object similar to a singular Hermitian metric $e^{-\vp _P}$ of $F$, defined as follows.  If $\xi$ is a local section of $F^*$ with no zeros, then
\[
\vp _P (x) := \log P (\xi _x,\xi _x).
\]
The reason this object is not necessarily a singular Hermitian metric is that we do not know $\log P (\xi _x,\xi _x)$ to be locally integrable.  Thus we cannot in general form the curvature current of this object.

\medskip

\noi {\bf Curvature.}
Given $P \in \sP _1 (X,F)$, we can define its G-curvatures (see \cite{dv}) as follows.  For each integer $N \ge 1$, the $N^{\text{th}}$ G-curvature of $P$ along $(v_1,...,v_N) \in (F^*)^N$ is the Hermitian matrix
\[
\Theta _N P (v_1,...,v_N) := \left (
\begin{array}{ccc}
P(v_1,\bar v_1) & \cdots & P(v_1,\bar v_N) \\
\cdot & \cdots & \cdot \\
\cdot & \cdots & \cdot \\
\cdot & \cdots & \cdot \\
P(v_N,\bar v_1) & \cdots & P(v_N,\bar v_N) \\
\end{array} \right )
\]
We let
\[
\sP _N (X,F)
\]
denote the set of Hermitian algebraic functions $P$ such that for any $v_1,...,v_N \in F^*$ the matrix
\[
\Theta _N P(v_1,...,v_N)
\]
is positive semi-definite.  Clearly
\[
\sP _{N+1}(X,F) \subseteq \sP _N (X,F).
\]
We write
\[
\sP _{\infty} (X,F) := \bigcap _{N \in \N} \sP _N (X,F).
\]
In \cite{dv} it is shown that
\begin{enumerate}
\item[(i)]$\sP _{\infty} (X,F)$ consists of all those Hermitian algebraic functions $P$ whose associated matrix $\sC _{\alpha \bar \beta}$ is positive semi-definite, and

\item[(ii)] there exists $N_0$ (depending on $F$) such that $\sP _{N_0} (X,F) = \sP _{\infty} (X,F)$.  (This conclusion uses the compactness of  $X$.)
\end{enumerate}

\noi In this paper we will only be concerned with the classes $\sP _1$, $\sP _2$ and $\sP _{\infty}$.

\medskip

\noi {\bf Sums of squared norms and quotients of squared norms.}
As we have just said, $\sP _{\infty} (X,F)$ consists of those Hermitian algebraic functions $P$ whose associated Hermitian matrix is positive semi-definite.  It follows that there is a basis $\{ s^1 ,...,s^k\}$ of $V_P$ such that
\[
P= \sum _{j=1} ^k |s^j|^2.
\]
For this reason, we refer to an element of $\sP _{\infty}$ as a {\it sum of squared norms}.

\begin{defn}\label{qos-defn}
A Hermitian algebraic function $P \in \sP _1(X,F)$ is said to be a quotient of squared norms if there exist sums of squared norms $P_1  \in \sP _{\infty} (X,F_1)$ and $P_2 \in \sP _{\infty} (X,F_2)$ such that
\[
F_1 = F\tensor F_2 \qquad {\rm and} \qquad P_1  = P \tensor P_2.
\]
The set of Hermitian algebraic functions associated to the line bundle $F \to X$ that are quotients of squared norms will be denoted $\cq (X,F)$.
\end{defn}

In the language of this section, we can phrase the following beautiful result of Catlin and D'Angelo \cite{cd3} as follows.

\begin{o-thm}\label{cd-qos-thm}
For any holomorphic line bundle $F \to X$,
\[
\sP _1 ^{\sharp} (X,F) \subset \cq (X,F).
\]
\end{o-thm}

\begin{rmk}
The special case of Theorem \ref{cd-qos-thm} corresponding to $X=\p _n$ was established independently by D. Quillen in \cite{q}.  Quillen's method of proof is on the surface very different from that of Catlin and D'Angelo.
\end{rmk}

Theorem \ref{cd-qos-thm} is a consequence of Theorem \ref{cd3-thm} given below.  To state the latter, we need the notion of Global Cauchy-Schwarz functions.

\medskip

\noi {\bf Global Cauchy-Schwarz function.}
The elements of the class $\sP _2$ are also called Global Cauchy-Schwartz functions, or simply {\it GCS-functions}, because the class $\sP _2$ is easily seen to consist precisely of those Hermitian algebraic functions $R$ that satisfy
\[
|R(v,\bar w)|^2 \le R(v,\bar v) R(w,\bar w).
\]
D'Angelo \cite{d1} proved that if $R \in \sP _2 (X,E)$ then the associated metric object $e^{-\vp _R}$ is a singular Hermitian metric for $E$ having semi-positive curvature current.  Calabi had also previously considered the functions $\vp _R$ in the paper \cite{cal}.

\begin{defn}\label{sgcs-defn}
Let $E \to X$ be a holomorphic line bundle.  A Hermitian algebraic function $R \in \sP_1(X,E)$ is said to be a strong GCS function (or simply SGCS function) if the following conditions hold.
\begin{enumerate}
\item[(S1)] $R \in \sP _1 ^{\sharp}(X,E) \cap \sP _2(X,E)$.
\item[(S2)] The matrix $\Theta _2 R(v,w)$ has determinant zero if and only if either $v=0$, $w=0$, or $v$ and $w$ lie in the same fiber of $E$.
\item[(S3)] The (smooth) metric $e^{-\vp _R}$ of $E$ associated to $R$ has strictly positive curvature.
\end{enumerate}
The class of SGCS metrics on $E$ will be denoted
$\sP ^S _2 (X,E).$
\end{defn}

\noi Definition \ref{sgcs-defn} is due essentially to Catlin and D'Angelo \cite{cd3}.  The main difference is that they do not include property (S3) in the definition, but make it a hypothesis in the following theorem.

\begin{o-thm}\label{cd3-thm}
\rm{(Catlin-D'Angelo \cite{cd3})}  Let $P \in \sP _1 ^{\sharp} (X,F)$ and $R\in \sP ^S _2 (X,E)$.  Then there is an integer $N_0$ such that for all $m \ge N_0$,
\[
R^m \tensor P \in \sP _{\infty} (X,E^m \tensor F).
\]
\end{o-thm}

\begin{rmk} By applying Theorem \ref{cd3-thm} with $P =1 \in \sP _1 ^{\sharp} (X,\co _X)$ we see that for all sufficiently large $m$, $R^m \in \sP _{\infty} (X,E^m)$, and thus Theorem \ref{cd-qos-thm} follows.
\end{rmk}

As mentioned in the Introduction, we shall give a proof of Theorem \ref{cd3-thm} in the last section of the paper.  

\section{Quotients of squared norms}\label{qos-section}

In this section we prove our main result, Theorem \ref{main}.  We start with the case where $X$ is one-dimensional (i.e., a compact Riemann surface).  Then we treat the case where the Hermitian algebraic function is free, and finally we establish the general case.

\medskip

\noi {\bf The Riemann surface case.}
Let $M$ be a compact Riemann surface and $\Lambda \to M$ a
holomorphic line bundle.

\begin{o-thm}\label{1d}
A Hermitian algebraic function $P \in \sP _1(M,\Lambda)$ is a quotient of squared norms if and only if there exist holomorphic line bundles $E_1,E_2 \to M$, a holomorphic section $\Phi \in H^0(M,E_1)$, and a positive Hermitian algebraic function $\tilde P$ on $M$ such that $\Lambda = E_1 \tensor E_2$ and
\begin{eqnarray}\label{split}
P = |\Phi |^2 \tilde P,& i.e.,& P(v,\bar w) = \left < \Phi ,v \right > \overline{\left < \Phi ,w\right >} \tilde P(v,\bar w).
\end{eqnarray}
\end{o-thm}

\begin{rmk} This theorem was proved by D'Angelo in the case $M = \p _1$ \cite{d2}.
Our proof is similar and slightly more efficient, but it is crucial for this paper that the surface $M$ is general.
\end{rmk}

\begin{proof}[Proof of Theorem \ref{1d}] The sufficiency of the
expression (\ref{split}) is a trivial consequence of Theorem
\ref{cd-qos-thm}.  To see the necessity, suppose $P$ is a quotient of squared norms.  Then
\begin{eqnarray}\label{decomp}
P\cdot \sum _{\alpha}  \left | s^{\alpha}\right |^2 = \sum
_{\beta}  \left | T^{\beta}\right |^2
\end{eqnarray}
for some holomorphic sections $\{s^{\alpha}\}$ and $\{
T^{\beta}\}$ of (different) line bundles.

Suppose now that $x \in M$.  Choose a local coordinate function $z$ near $x$ with $z(x)=0$, and let $\xi$ be a nowhere zero section of $\Lambda ^*$ near $x$.  Then by formula (\ref{decomp}) the function $p := P(\xi,\bar \xi)$ has the property that for some non-negative integer $m$, $p/|z|^{2m}$ is bounded and non-zero near $x$.  (If there exists $v \in \Lambda ^* _x - o_{\Lambda^*}(x)$ such that $P(v,\bar v) = 0$, then $m > 0$.)  This integer $m$ is unique and does not depend on the choice of local coordinate $z$.  We write
\[
m=:D_P(x).
\]
We associate to $P$ the divisor
\[
D_P= \sum _{x \in M} D_P(x)\cdot x.
\]
Evidently $D_P$ is independent of $\xi$ and $z$.  Let $L_{D_P}$ be the line bundle associated to $D_P$ and let $s_{D_P}$ be the canonical section of $L_{D_P}$ whose divisor is $D_P$.  Since $D_P$ is effective, i.e., $D_P(x) \ge 0$ for all $x \in M$, $s_{D_P}$ is holomorphic.

By construction,
\[
\tilde P = |s_{D_P}|^{-2} P \in \sP _1 (M, \Lambda \tensor L_{D_P} ^*)
\]
is a Hermitian algebraic function with no non-trivial zeros. This completes the proof.
\end{proof}

\ms

For the purpose of illustration, we prove Theorem \ref{main} for Riemann surfaces.

\ms

\begin{proof}[\underline{Necessity of basic zeros}]  Let $P \in \sP _1 (M,\Lambda)$.  First, suppose $P$ is a quotient of squared norms.
In view of Theorem \ref{1d}, $P = |s_D|^2 \tensor \tilde P$ for some positive Hermitian algebraic function $\tilde P$.  It follows that $s_D$ divides the support space $V_P$ of $P$.  Thus for any $s_D \tensor f \in V_P$ we have that
\[
\frac{|s_D\tensor f|^2}{P} = \frac{|f|^2}{\tilde P}
\]
is bounded.  That is to say, condition (\ref{bdd}) holds.

\ms

\noi \underline{\it Sufficiency of basic zeros.}  Conversely, suppose  condition (\ref{bdd}) holds for $P$.  If $V_P$ is free, then $P$ has no zeros and by Theorem \ref{cd-qos-thm} we are done.  Therefore, assume $V_P$ has a non-trivial base locus.  This base locus is a divisor $D$ on $M$ whose canonical holomorphic section $s_D$ divides every element of $V_P$.  It follows that $|s_D|^2$ divides $P$, and $\tilde P = P /|s_D|^2$ is a Hermitian algebraic function for the line bundle $E = \Lambda \tensor L_D^*$ such that $V_{\tilde P}$ is free.  Moreover,
\[
\frac{|s|^2}{P} = \frac{|s/s_D|^2}{\tilde P},
\]
And thus by condition (\ref{bdd}) $\tilde P$ has no zeros.
By Theorem \ref{cd-qos-thm} we are done.
\end{proof}

\begin{cor} \label{free-1d-rmk}
If $R$ is a Hermitian algebraic function on a Riemann surface and $V_R$ is free, then $R$ is a quotient of squared norms if and only if $R$ has no non-trivial zeros.
\end{cor}

\begin{ex}
Consider the Hermitian algebraic function $R \in \sP _1(\p _1, \h ^{\tensor 2})$ defined by 
\[
R(z,\bar w) = |w^2|^2 +(|z+w|^2-|z|^2)^2.
\]
Then the zero locus of $R$ is the point $w =0$ in $\p _1$.  However, the linear system $V_R = \text{span} _{\C} \{ w^2 , z^2 , (z+w)^2, z(z+w) \}$ is base point free, since it contains the base point free linear system $\text{span} _{\C} \{ z^2, w^2 \}$.  This show in particular that even if $\cz _R$ is a complex analytic set, $R$ need not be a quotient of squares.
\end{ex}

\medskip

\noi {\bf The free case.}  Corollary \ref{free-1d-rmk} generalizes to higher dimensions.

\begin{o-thm}\label{free-thm}
Let $P \in \sP _1 (X,F)$ be a Hermitian algebraic function such that $V_P$ is free. Then $P$ is a quotient of squared norms if and only if $P$ has no non-trivial zeros.
\end{o-thm}

\begin{proof}
The ``if'' direction follows from Theorem \ref{cd-qos-thm}.  We shall now prove the converse.  Suppose $P$ is a quotient of squared norms and $P(v,\bar v)=0$ for some $v\in F^*_x-o_{F^*}(x)$.  Since $X$
is projective algebraic, there is a Riemann surface $M \subset X$
passing through $x$. Let $\iota : M \emb X$ denote the natural
inclusion map.  Then the following facts are evident.
\begin{enumerate}
\item $\iota ^* V_P = V_{\iota ^* P}$, and thus $V_{\iota ^*P}$ is
free. \item $\iota ^* P$ is a quotient of squared norms. \item $\iota
^*P$ vanishes along $(\iota ^*F)_x$.
\end{enumerate}
But these three facts contradict Corollary \ref{free-1d-rmk}.  The proof is complete.
\end{proof}

\medskip

\noi {\bf An illustrative example.}
D'Angelo showed that the Hermitian algebraic function $P \in \sP _1 (\p _2 , \h ^{\tensor 4})$
given in homogeneous coordinates $[z_0,z_1,z_2]$ by
\[
P(z,\bar z) = |z_0|^8 + (|z_0z_2|^2 - |z_1|^4)^2
\]
is not a quotient of squared norms.  The two interesting aspects of this Hermitian algebraic function are that (i) its zero locus $\cz _P$ is the single point $[0,0,1]$, and thus $\cz _P$ is an analytic set, and (ii) every section in $V_P$ vanishes at $\cz _P$.

Nevertheless, D'Angelo showed that $P$ is not a quotient of squared norms.  His approach was to show that $P$ does not satisfy a certain necessary condition for being a quotient of squared norms: the jet pullback property.  We will recall the jet pullback property in Section \ref{jpp-section} below.  In the language of the present paper, D'Angelo passes the (rational) curve $\gamma : \p _1 \to \p _2$ given by
\[
\gamma  [x,y] = [x^2 , xy, xy + y^2]
\]
through $[0,0,1]$ and examines the restriction of $P$ to this curve.  He then deduces from Taylor expansions near $[0,0,1]$ that $P$ cannot be a quotient of squared norms.

The information about the influence of the Taylor expansion of the curve on $P$ is already contained in Theorem \ref{1d}.  To rephrase D'Angelo's proof in the language of the present paper, note that $\gamma ^* P \in \sP _1 (\p _1 , \h ^{\tensor 8})$ and
\begin{eqnarray*}
\gamma ^* P ((x,y),(\bar x,\bar y)) &=& |x|^{16}  + ( |x^3y + x^2y^2|^2 - |xy|^4)^2\\
&=& |x|^8 \left ( |x|^8 + ( |xy + y^2|^2 - |y|^4 )^2 \right )\\
&=& |x|^8 Q((x,y),(\bar x,\bar y)).
\end{eqnarray*}
One then sees that $Q$ is a free Hermitian algebraic function on $\p _1$ with a non-trivial zero at the point $[0,1] \in \p _1$.  By Corollary \ref{free-1d-rmk} $Q$ is not a quotient of squared norms.  Thus $\gamma ^* P$ is not a quotient of squared norms, and hence $P$ is not a quotient of squared norms.

\begin{rmk}
We leave it to the interested reader to show directly that condition \eqref{bdd} does not hold for $P$.  In fact, the set $\sB _P$ is empty.
\end{rmk}

In order to demonstrate the role of the resolution of singularities in the proof of Theorem \ref{main}, we shall use blowups to show that the Hermitian algebraic function $P$ is not a quotient of squared norms.  Since the zero locus of $P$ consists only of the point $[0,0,1] \in \p _2$, it is easiest to work in the affine chart $\{ z_2 \neq 0\}$.  Let $x_1 = z_0/z_2$ and $x_2 = z_1/z_2$.  Then
\[
P = |z_2|^8 (|x_1|^8 + (|x_1|^2 - |x_2|^4)^2).
\]
It suffices to work with the inhomogeneous polynomial
\[
p = |x_1|^8 + (|x_1|^2 - |x_2|^4)^2,
\]
since one can always recover the original Hermitian algebraic function by homogenization.  (For more on this point as well as other aspects of homogenization of polynomials with regard to Hermitian algebraic functions, see \cite{d2}.)

Now blow up the origin in $\C ^2$.  That is to say, let
\[
x_1 = y_1y_2 \quad \text{and} \quad x_2 = y_2.
\]
Then with $p_1$ denoting the blowup of $p$, one has
\[
p_1 = |y_2|^4\left ( |y_1^2 y_2|^4 + (|y_1|^2 - |y_2|^2)^2 \right ).
\]
Now take one more blowup, namely
\[
y_1 = t_1 \quad \text{and} \quad y_2 = t_1t_2.
\]
Then with $p_2$ denoting the blowup of $p_1$, one has
\[
p_2 = |t_1^2t_2|^4\left ( |t_1^2 t_2|^4 + (1- |t_2|^2)^2 \right ).
\]
We have shown that after two blowups, $p$ is transformed into a Hermitian algebraic function that is a product of a square (the term $|(t_1 ^2t_2)^2|^2$) and a free Hermitian algebraic function with a non-trivial zero, namely the point $(t_1,t_2) = (0,1)$ in the particular affine chart we are working with.  It follows from Theorem \ref{free-thm} that the Hermitian algebraic function $P$ could not have been a quotient of squared norms.

\medskip

\noi {\bf The general case.}
We remind the reader that we have assumed the absence of zeros in codimension 1.

Let $\mu : \tilde X \to X$ be a log resolution of $V_P$ as in Theorem \ref{log-res}.  Then there is an effective divisor $D$ in $\tilde X$, a free linear system $W \subset H^0(\tilde X,F\tensor L_D^*)$ and a section $s_D$ of the line bundle associated to $D$ such that for all $s \in \mu ^* V_P$,
\[
s= s_D \tensor t \text{ for some }t \in W.
\]
It follows that $\mu ^* P = |s_D|^2 \tensor \tilde P$ for some
Hermitian algebraic function $\tilde P \in \sP _1 (X,\mu ^* F \tensor L_D^*)$ such that
$V_{\tilde P}$ is free.

We now show the necessity of condition (\ref{bdd}).  To this end,
if $P$ is a quotient of squared norms, then so is $\tilde P$.  Since
$V_{\tilde P}$ is free, by Theorem \ref{free-thm} we see that
$\tilde P$ has no zeros.  It follows that if $s \in V_P$ and $\mu
^* s = s_D \tensor t$ for some $t \in V_{\tilde P}$, then the
quotient $|t|^2/ \tilde P$ is bounded. But then the quotient
\[
\frac{|\mu ^* s|^2}{\mu ^* P}= \frac{|s_D|^2 |t|^2}{|s_D|^2 \tilde
P}
\]
is bounded on $\tilde X$.  Since $\mu$ is surjective, we see that condition (\ref{bdd}) must hold.

Next we turn to the sufficiency of condition (\ref{bdd}).  To this end, let $E \to X$ be an ample line bundle.  Then $\mu ^*E \to \tilde X$ is big.  Replacing $E$ be a large power if needed, we may assume that $\mu ^*E = A \tensor L_Z$ where $A \to \tilde X$ is a very ample line bundle and $L_Z$ is the line bundle associated to an effective divisor $Z$.  Let $s_Z$ be the canonical holomorphic section of $L_Z$ with zero divisor $Z$, and let $\hat s^1,...,\hat s ^N$ be global sections of $A \to \tilde X$ such that the Hermitian algebraic function
\[
\hat R := |\hat s^1|^2 + ... +|\hat s^N|^2
\]
is SGCS.  (The condition of being in $\sP _2$ is trivial in this case, and the positivity of the curvature follows from the very-ampleness of $A$ and the positivity of the Fubini-Study metric.)  By construction, there is a Hermitian algebraic function $R \in \sP _1(X,E)$ such that
\[
\mu ^* R = |s_Z|^2 \tensor \hat R.
\]
Now,
\[
\mu ^* (R^m \tensor P) = |s_Z ^{\tensor m} \tensor s_P|^2 \tensor \hat R^m \tensor \tilde P.
\]
It follows from Theorem \ref{cd3-thm} that $\mu ^* (R^m \tensor P) $ is a quotient of squared norms, and thus the same is true of $R^m \tensor P$.  The proof of Theorem \ref{main} is complete.  \qed

\begin{rmk}
Note that the Hermitian algebraic function $R$ might not be positive if $P$ is not positive.  In fact, it is not always possible to choose $R$ positive when $P$ has non-trivial zeros.
\end{rmk}

\section{The jet pullback property}\label{jpp-section}

Let us recall the definition of D'Angelo's jet pullback property.

\begin{defn}  A Hermitian algebraic function $P \in \sP _1 (X,F)$ is said to have the {\it jet pullback property} if for every compact Riemann surface $M$ and holomorphic map $h : M \to X$ the following holds.  If $z$ is a local coordinate in a neighborhood of $x \in M$ and $\xi$ is a nowhere zero holomorphic section of $F$ in a neighborhood of $h(x)$, then the function $h^* (P(\xi,\xi))$ has Taylor expansion
\[
h^* (P(\xi,\xi)) (z) = c |z|^{2d} + \text{higher order terms}
\]
for some $c >0$.
\end{defn}

\begin{rmk}
Strictly speaking, in D'Angelo's definition of the jet pullback property  one is only allowed to use rational curves, i.e., $M = \p ^1$.  We shall take this stronger definition, since rational curves might not be so plentiful on $X$.  Indeed, it may happen that $X$ has no rational curves at all.
\end{rmk}

It is easy to see that if $P$ is a quotient of squared norms then $P$ has the jet pullback property.  In \cite{d-carrus} it is asked whether the converse is true.  In this section we answer this question in the affirmative.

\begin{prop}\label{jpp-qos}
If $P \in \sP _1(X,F)$ has the jet pullback property then $P$ is a quotient of squared norms.
\end{prop}

\begin{proof}
Suppose $P$ has the jet pullback property.  We shall show that $P$ has only basic zeros, and then use Theorem \ref{main}.

\noi (i) As in the previous section, let $\mu : \tilde X \to X$ be a resolution of singularities such that $\mu ^* P = |s|^2 \tilde P$ for some holomorphic section $s$ and free Hermitian algebraic function $\tilde P$.

\ms

\noi (ii) Observe that $\mu ^* P$ has the jet pullback property.  Indeed, if $h : M \to \tilde X$ is a curve, then
\[
h ^* (\mu ^* P) = (\mu \circ h)^* P,
\]
and thus it is clear from the jet pullback property for $P$ that the lowest order term of the Taylor expansion for $h^* (\mu ^* P)$ at any point of $M$ has the right form.

\ms

\noi (iii) Notice that since $|s|^2 \tilde P$ has the jet pullback property, so does $\tilde P$.

\ms

\noi (iv)  If the free Hermitian algebraic function $\tilde P$ has the jet pullback property, then $\tilde P$ must be a positive Hermitian algebraic function.  Indeed, since $\tilde P$ is free, the pullback $h^* \tilde P$ of $\tilde P$ to any curve $h : M \to \tilde X$ is free.  Moreover, it is immediate that $h^* \tilde P$ has the jet pullback property.  But since $M$ is a curve, the jet pullback property for $h^* \tilde P$ means that the zero set of $h^* \tilde P$ is a divisor.   If $\tilde P$ is not a positive Hermitian algebraic function, then by choosing the curve $h :M \to X$ such that $h(M)$ passes through a zero of $\tilde P$, we see that the zero set of $h^* \tilde P$ is non-empty.  It follows as in the proof of Theorem \ref{1d} that every section in the support space of $h^* \tilde P$ is divisible by a global section on $M$, and this contradicts that $h^* \tilde P$ is free.

\noi (v)  If $\tilde P$ is a positive Hermitian algebraic function, then clearly
\[
\mu ^*\left ( \frac{\gmod{P}}{P} \right ) = \frac{\gmod{\tilde P}}{\tilde P}
\]
is uniformly bounded on $\tilde X$.  Since $\mu$ is surjective, we see that $\frac{\gmod{P}}{P}$ is uniformly bounded on $X$.  By Theorem \ref{main} we are done.
\end{proof}

Using Proposition \ref{jpp-qos} it is easy to see the following, alternate characterization of quotients of squared norms.

\begin{prop}
A Hermitian algebraic function $P \in \sP _1(X,F)$ is a quotient of squared norms if and only if for every holomorphic map $h : M \to X$ of a compact Riemann surface $M$ into $X$, $h^*P \in \sP _1 (M, f^*F)$ is a quotient of squared norms.
\end{prop} 

\begin{proof}
One direction is immediate: if $P$ is a quotient of squared norms then the pullback of $P$, by a holomorphic map, to any complex manifold is also a quotient of squared norms.  Conversely, suppose $P$ is not a quotient of squared norms.  Then by Proposition \ref{jpp-qos} there is a Riemann surface $M$ and a holomorphic map $h : M \to X$ such that $P_o := h^*P \in \sP _1 (M, h^* F)$ does not have the jet pullback property.   Thus again by Proposition \ref{jpp-qos} $P_o$ is not a quotient of squared norms.  The proof is complete.
\end{proof}

\section{Bergman kernels.  The Proof of Theorem \ref{cd3-thm}}

\noi In this final section we give an alternate proof of the Catlin-D'Angelo Theorem \ref{cd3-thm}.

\medskip

\noi {\bf Spaces of $L^2$-sections and the generalized Bergman kernel.}
Let $(X, \omega)$ be a K\" ahler manifold and let 
\[
{\Omega} := \frac{\omega ^n}{n!(2\pi)^n}
\] 
be the associated volume form.  Suppose given a holomorphic line bundle $E \to X$ with singular Hermitian metric $e^{-\vp}$.  For smooth sections $f,g$ of $E$, we define
\[
(f,g) := \int _X f\bar g e^{-\vp} \Omega.
\]
With this inner product, we denote by $L^2 (e^{-\vp})$ the Hilbert space completion of the set of smooth sections $f$ of $E \to X$ such that
\[
||f||^2 = \int _X |f|^2 e^{-\vp} \Omega < +\infty.
\]

The Hilbert space $L^2 (e^{-\vp})$ has a closed, finite-dimensional subspace
\[
\sa ^2 (e^{-\vp}) := L^2 (e^{-\vp}) \cap H^0(X,E).
\]
The projection $\K _{\vp} :L^2 (e^{-\vp}) \to \sa ^2 (e^{-\vp})$
is an integral operator, called the (generalized) Bergman
projection.  Its kernel, $K_{\vp}$, is given by the formula
\[
K_{\vp}(x,y) = \sum _{\alpha =1} ^N s^{\alpha}(x) \tensor
\overline{s^{\alpha}(y)},
\]
where $\{ s^1,...,s^N \}$ is an orthonormal basis for $\sa ^2
(e^{-\vp})$.  

\medskip

\noi {\bf The Bergman kernel of $H^0(X,E^m \tensor F)$.}
Let $E \to X$, $F \to X$ be holomorphic line bundles, and fix
Hermitian algebraic functions $R \in \sP _2 ^S (X,E)$ and $P \in \sP_1
^{\sharp}(X,F)$.  We write $e^{-\vp}$ for the smooth metric of $E$ associated to $R$ and $e^{-\psi}$ for the smooth metric of $F$ associated to $P$.  In a similar fashion, given local non-vanishing holomorphic sections $\xi , \eta$ of $E^* \to X$ near points $x_o,y_o \in X$ respectively, we write $r(x,y) := R(\xi _x,\bar \eta _y)$, and likewise $p(x,y)$ for the local expression of $P$.  Thus
\[
\vp (x) := \log r(x,x) \quad \text{and} \quad \psi (x) := \log
p(x,x)
\]
provided we use $\xi = \eta$ when $x$ and $y$ are close.

From here on we fix our K\"ahler form to be 
\[
\omega := \ii \di \dbar \vp.
\]
We have the Hilbert spaces $L^2 (e^{-(m\vp + \psi)})$
and $\sa ^2 (e^{-(m\vp + \psi)})$ from the previous paragraph.  Because $R$ and $P$ are positive Hermitian algebraic functions,
$\sa ^2 (e^{-(m\vp + \psi)}) = H^0(X,E^m\tensor F),$ as vector spaces, i.e., every global holomorphic section is $L^2$.  We write
\[
K^{(m)} := K_{m \vp + \psi},
\]
for the Bergman kernel of $\sa ^2(e^{-(m\vp + \psi)})$.

The key result needed is the following theorem.

\begin{o-thm} \label{Bergman-asymp}
There exist smooth functions $b_j(x,y)$, $j=1,2,...n$, such that 
\[
K^{(m)} (x,y) -  m^n r(x,y)^m p(x,y)\left ( 1-\sum _{j=1} ^{n} \tfrac{1}{m^j} b_j(x,\bar y) \right ) = O(\tfrac{1}{m}).
\]
\end{o-thm}

In the case where $p \equiv 1$, Theorem \ref{Bergman-asymp} is Theorem 2 in \cite{c} or (the polarization of) Theorem 1 in \cite{z}.  
As mentioned in the introduction, the proofs of both of these theorems make use of the Boutet de Monvel-Sj\" ostrand Theorem \cite{bs}.  A more direct proof using more elementary methods has recently appeared in \cite{bbs}.  In Section 2.5 of that paper, it is indicated  how to modify the argument to the case of general $p$ stated above.  In \cite{bbs} only the local version of Theorem \ref{Bergman-asymp} is proved (as Theorem 3.1) , but a partition of unity argument easily extends their result to the one we have stated here.

\medskip

\noi {\bf Proof of Theorem \ref{cd3-thm}.}
Choose orthonormal bases $s^1,...,s^{N_m}$ for the spaces $L^2(e^{-(m\vp +\psi)})$, and let 
\[
s = \sum _{i=1} ^{N_m} a^i s^i
\]
be a unit vector.  If we denote by $\sC ^{(m)} _{\alpha \bar \beta}$ the matrix associated to the Hermitian algebraic function $R^m P$, then we have 
\[
\sC ^{(m)} _{\alpha \bar \beta} a^{\alpha} \bar a^{\beta} = \int _X \int _X \frac{r(x,y)^m p(x,y) s(y) \overline{s(x)}}{r(x,x)^m p(x,x) r(y,y)^mp(y,y)} \Omega (y) \Omega (x), 
\]
while 
\begin{equation*}\label{reprod}
\delta _{\alpha \bar \beta} a^{\alpha} \bar a^{\beta} = \int _X \int _X \frac{k^{(m)}(x,y)  s(y) \overline{s(x)}}{r(x,x)^m p(x,x) r(y,y)^mp(y,y)} \Omega (y) \Omega (x).
\end{equation*}

Now consider the linear operators $L_i ^{(m)} : \C ^{N_m} \to \C ^{N_m}$ defined by 
\[
(L_i ^{(m)}) _{\alpha \bar \beta} a^{\alpha} \bar a^{\beta} = \int _X \int _X \frac{r(x,y)^m p(x,y) b_i(x,\bar y) s(y) \overline{s(x)}}{r(x,x)^m p(x,x) r(y,y)^mp(y,y)} \Omega (y) \Omega (x), \quad i=1,..., n.
\]
These operators are clearly comparable to $\sC ^{(m)}$ uniformly in $m$, in the sense that there is a constant $C > 0$, independent of $m$ such that for all $1 \le i \le n$ and all $a \in \C ^{N_m}$, 
\[
-C \left < L_i ^{(m)} a, a\right > \le \left < \sC ^{(m)} a, a\right > \le C \left < L_i ^{(m)} a, a\right >.
\]
Applying Theorem \ref{Bergman-asymp}, we see that 
\[
||a||^2 = m^n \left < \sC ^{(m)} a, a\right > - \sum _{i=1} ^{n} m^{n-i} \left < L_i ^{(m)} a, a\right > + O(\frac{1}{m}).
\]
It follows that there is a polynomial $P(m)$ of degree $n$,  with leading coefficient $1$, such that 
\[
||a||^2 \le P(m)  \left < \sC ^{(m)} a, a\right >.
\]
Let $N_o := \min \{ m_o ; P(m) > 0 \text{ for all }m \ge m_o \}.$  Then for $m \ge N_o$, we have 
\[
\left < \sC ^{(m)} a, a\right > > 0.
\]
This completes the proof of Theorem \ref{cd3-thm}.


\begin{thebibliography}{99}
\bibitem[AD]{ad} Abramowitz, D., De Jong, A. J.,
{\it Smoothness, semistability, and toroidal geometry. }
J. Algebraic Geom. 6 (1997), no. 4, 789--801.

\bibitem[BBS]{bbs} Berman, R., Berndtsson, B., Sj\" ostrand, J., 
{\it Asymptotics of Bergman kernels.}  Preprint, 2005.  CV/0506367 v1.

\bibitem[BP]{bp} Bogomolov, F., Pantev, T.,
{\it Weak Hironaka theorem.}
Math. Res. Lett.  3  (1996),  no. 3, 299--307.

\bibitem[BS]{bs} Boutet de Monvel, Sj\"ostrand, J.,
{\it Sur la singularit\' e des noyaux de Bergman et de Szeg\" o.} (French)  Journ\' ees: \' Equations aux D\' eriv\' ees Partielles de Rennes (1975),  pp. 123--164. Asterisque, No. 34-35, Soc. Math. France, Paris, 1976.

\bibitem[Cal]{cal} Calabi, E., {\it Isometric imbedding of complex manifolds. }  Ann. of Math. (2) 58, (1953). 1--23. 

\bibitem[C]{c} Catlin, D.,
{\it The Bergman kernel and a theorem of Tian.}
Analysis and geometry in several complex variables (Katata, 1997),
1--23, Trends Math., Birkh\" auser Boston, Boston, MA, 1999.

\bibitem[CD1]{cd1} Catlin, D., D'Angelo, J.,
{\it A stabilization theorem for Hermitian forms and applications to holomorphic mappings.}
Math. Res. Lett. 3  (1996),  no. 2, 149--166.

\bibitem[CD2]{cd3} Catlin, D., D'Angelo, J.,
{\it An isometric embedding theorem for holomorphic bundles.} Math. Res. Lett. 6 (1999) no. 1, 43--60.

\bibitem[D1]{d-carrus} D'Angelo, J.,
{\it Inequalities from complex analysis.}
Carus Monograph 28. MAA,  2002.

\bibitem[D2]{d1} D'Angelo, J.,
{\it Bordered complex Hessians.}
J. Geom. Anal.  11  (2001),  no. 4, 561--571.

\bibitem[D3]{d2} D'Angelo, J.,
{\it Complex Variables Analogues of Hilbert's Seventeenth Problem}
To appear in Internat. J. Math.

\bibitem[DV]{dv} D'Angelo, J., Varolin, D.,
{\it Positivity conditions for Hermitian symmetric functions.}
Asian J. Math, 7, No. 4 (2003) 1-18.

\bibitem[F]{f} Fefferman, C.,
{\it The Bergman kernel and biholomorphic mappings of pseudoconvex domains.}
Invent. Math. 26 (1974), 1--65.

\bibitem[H]{h} Hironaka, H.,
{\it Resolution of singularities of an algebraic variety over a field of characteristic zero. I, II.}
Ann. of Math. (2) 79 (1964), 109-203, 205-326.

\bibitem[L]{rlaz} Lazarsfeld, R.,
{\it Positivity in Algebraic Geometry I and II}
Springer 2004.

\bibitem[P]{p} Paranjape, K.,
{\it  The Bogomolov--Pantev resolution, an expository account.} New Trends in Algebraic Geometry (ed. K. Hulek, F. Catanese, C. Peters, M. Reid), pages 347-358, Cambridge University Press 1999.   (See also  math.AG/9806084)


\bibitem[Q]{q} Quillen, D.,
{\it On the representation of Hermitian forms as sums of squares.}  Invent. Math. 5 (1968) 237--242.

\bibitem[Z]{z} Zelditch, S.,
{\it Szeg\" o kernels and a theorem of Tian.}
Internat. Math. Res. Notices  1998,  no. 6, 317--331.


\end{thebibliography}
\end{document}